\documentclass[11pt]{article}
\usepackage{mathrsfs}
\usepackage{amssymb}
\usepackage{amsmath}
\usepackage[all]{xy}
\def\Hom{{\rm Hom}}
\def\Ext{{\rm Ext}}
\def\Add{{\rm Add}}

\def\Im{\mathop{\rm Im}\nolimits}
\def\Ker{\mathop{\rm Ker}\nolimits}
\def\Coker{\mathop{\rm Coker}\nolimits}

\def\fd{\mathop{\rm fd}\nolimits}
\def\id{\mathop{\rm id}\nolimits}
\def\pd{\mathop{\rm pd}\nolimits}
\def\gldim{\mathop{\rm gl.dim}\nolimits}
\def\PD{{\rm PD}}
\def\ID{{\rm ID}}
\def\P{{\rm P}}
\def\I{{\rm I}}
\def\G{{\rm G}}
\def\w{{\rm w}}

\def\inf{\mathop{\rm inf}\nolimits}
\def\sup{\mathop{\rm sup}\nolimits}

\def\-{{\text{-}}}

\title{\Large \bf Global dimensions of rings with respect to a semidualizing module
\thanks{2000 Mathematics Subject Classification: 13D02, 13D05, 13D07, 18G25.}
\thanks{Keywords: semidualizing modules; Strongly $G_C$-projective and injective
modules; $G_C$-projective and injective dimensions; $G_C$-(weak) global dimension of rings.
}}
\author{Guoqiang Zhao$^{a}$\thanks{Corresponding author.}, Juxiang Sun$^b$\\
{\footnotesize$^a$\it School of Science, Hangzhou Dianzi University,
Hangzhou, 310018, China}\\{\footnotesize$^b$\it School of Mathematics and
Information Science, Shangqiu Normal University, Shangqiu, 476000, China}\\
{\footnotesize\it E-mail address: gqzhao@hdu.edu.cn, sunjx8078@163.com}}
\date{ }
\begin{document}
\baselineskip=18pt \maketitle
\begin{abstract}
In this paper, the notion of strongly $G_C$-projective and injective
modules is introduced, where $C$ is a semidualizing module. Using these modules
we can obtain a new characterization of $G_C$-projective and injective modules,
similar to the one of projective modules by the free modules.
We then define and study the global dimensions of rings relative to a semidualizing
module $C$, and prove that the global $G_C$-projective dimension of a ring $R$ is
equal to the global $G_C$-injective dimension of $R$.

\end{abstract}

\vspace{0.5cm}

\centerline{\large \bf 1. Introduction}

\vspace{0.2cm}


Over a Noetherian ring, Foxby [5] introduced the notion of semidualizing modules,
which provided a common generalization of a dualizing module and a
free module of rank one. Golod [6] and Vasconcelos [13] furthered the study of
semidualizing modules. By using these modules, Golod defined
the $G_C$-dimension, a refinement
of projective dimension, for finitely generated modules. When $C = R$,
this recovers the $G$-dimension introduced by Auslander and Bridger in [1].
Motivated by Enochs and Jenda$^{'}$s extensions
in [4] of $G$-dimension, Holm and J$\phi$gensen [8]
have extended the $G_C$-dimension to arbitrary modules over a Noetherian ring
(where they used the name of $C$-Gorenstein projective dimension). This also enables
them to give the dual notion. Then White [14] extended these concepts to the non-Noetherian setting,
named $G_C$-projective and $G_C$-injective dimension,
and showed that they share many common properties with
the Gorenstein homological dimensions extensively studied in recent decades.


It is well-known that, the classical global dimensions of rings play
an important role in the theory of rings. Recently, Bennis and Mahadou [3]
defined the global Gorenstein projective dimension of a ring
$R$ and the global Gorenstein injective dimension of $R$, and proved that they are equal for
any associative ring by using the properties of strongly Gorenstein projective and injective
modules. For more details of these two modules, see [2].
Based on the above results, in this paper we mainly study the global dimensions
of a ring $R$ with respect to a semidualizing module.

This paper is organized as follows.

In Section 2, we give some definitions and some preliminary results.

In Section 3, we introduce and study the
strongly $G_C$-projective and injective modules.
The main result of this section is that an $R$-module is $G_C$-projective
(resp. injective) if and only if it is a direct summand of a strongly $G_C$-projective
(resp. injective) $R$-module. We then give some equivalent characterizations of
the strongly $G_C$-projective and injective modules, and show that the class of
strongly $G_C$-projective (resp. injective) modules is between the class of projective,
$C$-projective (resp. injective) modules and the class of $G_C$-projective (resp. injective) modules.

In Section 4, relative to a semidualizing module $C$, we firstly define the
$\mathcal{P}_{C}$-projective and $\mathcal{I}_{C}$-injective dimension of a ring $R$,
and prove that they are both equal to the classical global dimension of $R$. Next, we
define and investigate the $G_{C}$-projective and $G_{C}$-injective dimension of
a ring $R$. The main result of this
paper is that the $G_{C}$-projective and $G_{C}$-injective dimension of a ring $R$
coincide, and we call the common value the $C$-Gorenstein
global dimension of $R$. Then we discuss the relations
between the $C$-Gorenstein global dimension of a ring $R$ with
other global dimensions of $R$. At the end of this section,
we study the behavior of modules over rings of finite
$C$-Gorenstein global dimension, and give a partial answer to the
question posed by Takahashi and White in [12].

\vspace{0.5cm}

\centerline{\large\bf 2. Preliminaries}

\vspace{0.2cm}

Throughout this work $R$ is a commutative ring with unity. For an $R$-module $M$,
we use id$_R(M)$, pd$_{R}(M)$ and fd$_{R}(M)$
to denote the injective dimension, projective dimension and flat dimension of
$M$, respectively. We use $\gldim(R)$ to denote the classical global dimension of $R$.

Semidualizing modules, defined next, form the basis for our categories of interest.

{\bf Definition 2.1} ([14]) An $R$-module $C$ is {\it semidualizing} if

(a) $C$ admits a degreewise finite projective resolution,

(b) The natural homothety map $R\rightarrow\Hom_R(C,C)$ is an isomorphism, and

(c) Ext$^i_R (C,C) = 0$ for any $i\geq 1$.

From now on, $C$ is a semidualizing module.

{\bf Definition 2.2} ([9]) An $R$-module is called {\it
$C$-projective} if it has the form $C\otimes_R P$ for some
projective $R$-module $P$. An $R$-module is called {\it
$C$-injective} if it has the form $\Hom_R(C, I)$ for some injective
$R$-module $I$. Set
$$\mathcal{P}_C=\mathcal{P}_C(R) = \{C\otimes_R P |_RP\ is\ projective\}, and$$
$$\mathcal{I}_C=\mathcal{I}_C(R) = \{\Hom_R(C, I) |_RI\ is\ injective\}.$$

Let $\mathcal{C}$ be a subclass of $R$-modules. Recall that a sequence
of $R$-modules $\mathbf{L}$ is called $\Hom_R(-, \mathcal{C})$
(resp. $\Hom_R(\mathcal{C}, -)$) exact if the sequence $\Hom_R(\mathbf{L}, C^\prime )$
(resp. $\Hom_R(C^\prime, \mathbf{L})$) is exact for any $C^\prime\in \mathcal{C}$.

{\bf Definition 2.3} ([14]) A complete $\mathcal{PP}_C$-resolution is a
$\Hom_R(-,\mathcal{P}_C)$ exact exact sequence of $R$-modules
$$\mathbf{X}= \cdots\rightarrow P_1\rightarrow
P_0\rightarrow C\otimes_RP^0\rightarrow C\otimes_RP^1\rightarrow \cdots$$
with $P_i$ and $P^i$ are projective $R$-modules. An $R$-module $M$ is
called {\it $G_C$-projective} if there exists
a complete $\mathcal{PP}_C$-resolution as above with
$M\cong \Coker(P_1\rightarrow P_0)$. Set
$$\mathcal{GP}_C(R) = \mathrm{the\ class\ of\ G_C\-projective\ R\-modules}.$$

A complete $\mathcal{I}_C\mathcal{I}$-resolution and $G_C$-injective module are
defined dually.

From [14, proposition 2.6], we know that every projective and
$C$-projective $R$-module are $G_C$-projective. Now, we give a "non-trivial" example
of $G_C$-projective $R$-module.

{\bf Example 2.4} Assume $R$ is a Gorenstein artin algebra with $\gldim(R)=\infty$. Let $C = \oplus I_j$,
where $I_j$ are all the indecomposable and non-isomorphic direct summands of modules
appeared in the minimal injective resolution of $R$. Then $C$ is a semidualizing module.
In this case, every finitely generated $R$-module is $G_C$-projective. While
the class of finitely generated $C$-projective $R$-modules is just the class of all
finitely generated injective $R$-modules. However, it is clear that there exists an
$R$-module which is not projective and injective.

Let $\mathcal{G}^2\mathcal{P}_C(R)$ = $\{$ $A$ is an $R$-module $|$ there exists a $\Hom_R(-,\mathcal{P}_C)$ exact exact
sequence of $R$-modules $\cdots\rightarrow G_1\rightarrow
G_0\rightarrow G^0\rightarrow G^1\rightarrow\cdots$ with all $G_i$ and $G^i$ in $\mathcal{G}\mathcal{P}_C(R)$ and $A\cong$ Im$(G_0\rightarrow G^0)\}$.

The following result means that an iteration of the procedure used to define the $G_C$-projective
modules yields exactly the $G_C$-projective modules.

{\bf Lemma 2.5} ([10, Theorem 2.9]) $\mathcal{G}^2\mathcal{P}_C(R)$ = $\mathcal{GP}_C(R)$.

{\bf Definition 2.6} ([14]) Let $\mathcal{X}$ be a class of $R$-modules and
$M$ an $R$-module. An {\it $\mathcal{X}$-resolution} of $M$ is an
exact sequence of $R$-modules as follows:
$$\cdots\rightarrow
X_n\rightarrow\cdots\rightarrow X_1\rightarrow X_0\rightarrow
M\rightarrow 0$$ with each $X_i$ $\in$ $\mathcal{X}$ for any $i\geq
0$. The $\mathcal{X}$-${\it projective\ dimension}$ of $M$ is the
quantity
$$\mathcal{X}\-\pd_{R}(M) = \inf\{\sup \{n\geq0 | X_n\neq0\} | X is\ an\ \mathcal{X}\-resolution\ of\ M\}.$$
The {\it $\mathcal{X}$-coresolution} and $\mathcal{X}$-${\it
injective\ dimension}$ of $M$ are defined dually. We write
$G_{C}\-\pd_R(M)$ = $\mathcal{GP}_C(R)\-\pd_R(M)$.

{\bf Lemma 2.7} {\it If $\sup\{G_{C}\-\pd_R(M)| M\ is\ an\
R\-module\}<\infty$, then, for an integer $n$, the following
are equivalent:

(1) $\sup\{G_{C}\-\pd_R(M)| M\ is\ an\ R\-module\} \leq n$,

(2) $\id_R(N) \leq n$ for every $R$-module $N$ with finite
$\mathcal{P}_C$-projective dimension.}

{\it Proof.} Use [14, Proposition 2.12] and [11, Theorem 9.8].$\hfill{\square}$

{\bf Definition 2.8} ([8]) $M$ is called {\it $G_C$-flat} if
there is an exact sequence of $R$-modules
$$\mathbf{X}= \cdots\rightarrow F_1\rightarrow
F_0\rightarrow C\otimes_RF^0\rightarrow C\otimes_RF^1\rightarrow \cdots$$
with $F_i$ and $F^i$ are flat $R$-modules, such that $M\cong \Coker(F_1\rightarrow F_0)$
and $\Hom_R(C, I)\otimes_R\mathbf{X}$ is still exact for any injective $R$-module $I$.
We define $G_{C}\-\fd_R(M)$ analogously to $G_{C}\-\pd_R(M)$.

 \vspace{0.5cm}

\centerline{\large\bf 3. Strongly $G_C$-projective and injective modules}

\vspace{0.2cm}

We denote $\Add_R M$ the subclass of $R$-modules consisting of all
modules isomorphic to direct summands of direct sums of copies of
$M$. By [10, Proposition 2.4], $\mathcal{P}_C = \Add_R C$.

{\bf Definition 3.1}  An $R$-module $M$ is called {\it strongly
$G_C$-projective}, if there exists an exact sequence of $R$-modules
$$\mathbf{D}= \cdots\stackrel{f}{\longrightarrow} D\stackrel{f}{\longrightarrow}
D\stackrel{f}{\longrightarrow} D\stackrel{f}{\longrightarrow} \cdots$$
with $D\in \Add_R ( C\oplus R )$, such that $M\cong \Ker\, f$ and
$\Hom_R( \mathbf{D}, \mathcal{P}_{C})$ is still exact.

When $C = R$, it is just the strongly Gorenstein projective module
introduced in [2].

Strongly $G_C$-injective modules are defined dually. In the following,
we only deal with the strongly $G_C$-projectivity of modules. The results about
strongly $G_C$-injective modules have a dual version, and we omit them.

From definition, we immediately have:

{\bf Proposition 3.2} {\it The class of strongly $G_C$-projective modules is closed under direct sums.}

The principal role of these modules is to give a simple characterization of the
$G_C$-projective modules, as follows:

{\bf Theorem 3.3} {\it An $R$-module is $G_C$-projective if and
only if it is a direct summand of a strongly $G_C$-projective $R$-module.}

{\it Proof.} $(\Leftarrow)$ From Lemma 2.5 and [14, Proposition 2.6],
it is easy to see that every strongly $G_C$-projective module is $G_C$-projective.
Since the class of $G_C$-projective modules is closed under
direct summands by [14, Theorem 2.8], the assertion follows immediately.

$(\Rightarrow)$ Let $M$ be a $G_C$-projective $R$-module. Then the definition
gives rise to a $\Hom_R ( -, \mathcal{P}_{C})$ exact exact sequence of $R$-modules:
$$X = \cdots\longrightarrow C_1\stackrel{d_1}{\longrightarrow} C_0\stackrel{d_0}
{\longrightarrow} C_{-1}\stackrel{d_{-1}}{\longrightarrow}
C_{-2}\longrightarrow\cdots$$ with all $C_i\in \Add R\cup \Add C$
and $M\cong \Im(C_0\rightarrow C_{-1})$.

For each $n\in \mathbb{Z}$, let $\Sigma^nX$ be the exact complex
obtained from $X$ by increasing all index by $n$: $(\Sigma^nX)_i =
X_{i-n}$ and $d^{\Sigma^nX}_i = d_{i-n}$ for all $i\in \mathbb{Z}$.

Then we obtain an exact complex
$$\bigoplus\Sigma^nX = \cdots\longrightarrow \oplus C_i\stackrel{\oplus d_i}{\longrightarrow} \oplus C_i\stackrel{\oplus d_i}
{\longrightarrow} \oplus C_i\stackrel{\oplus d_i}{\longrightarrow} \oplus C_i\longrightarrow\cdots$$


Clearly, $\oplus C_i\in\Add_R ( C\oplus R )$. Since $\Hom_R ( \bigoplus\Sigma^nX, \mathcal{P}_{C})\cong\prod\Hom_R (\Sigma^nX, \mathcal{P}_{C})$, the complex $\bigoplus\Sigma^nX$ is also
$\Hom_R ( -, \mathcal{P}_{C})$ exact. Thus
$M$ is a direct summand of the strongly $G_C$-projective module
$\Im(\oplus d_i)$, as desired.$\hfill{\square}$

The next result gives a simple characterization of the strongly
$G_C$-projective modules.

{\bf Proposition 3.4} {\it For any $R$-module $M$, the following are equivalent:

(1)$M$ is  strongly $G_C$-projective;

(2)There exists a short exact sequence $0\rightarrow M\rightarrow
D\rightarrow M\rightarrow 0$, with $D\in \Add_R ( C\oplus R )$, and
$\Ext^{1}_R(M, J)=0$ for any $R$-module $J$ with finite
$\mathcal{P}_C$-projective dimension ( or for any $C$-projective $R$-module $J$ );

(3)There exists a short exact sequence $0\rightarrow M\rightarrow
D\rightarrow M\rightarrow 0$, with $D\in \Add_R ( C\oplus R )$, and
$\Ext^{i}_R(M, J)=0$ for some integer $i > 0$ and for any $R$-module $J$
with finite $\mathcal{P}_C$-projective dimension.}

{\it Proof.} (1) $\Rightarrow$ (2) Follows from definition and
[14, Proposition 2.12], and (2) $\Rightarrow$ (3) is trivial.

(3) $\Rightarrow$ (1) Let $0\rightarrow M\rightarrow
D\rightarrow M\rightarrow 0$ be the short exact sequence
with $D\in \Add_R ( C\oplus R )$. Then, for for any $R$-module $J$
with finite $\mathcal{P}_C$-projective dimension and all $j > 0$, we
have the long exact sequence:
$$0 = \Ext^{j}_R(D, J)\rightarrow \Ext^{j}_R(M, J)\rightarrow
\Ext^{j+1}_R(M, J)\rightarrow \Ext^{j+1}_R(D, J) = 0$$
Then $\Ext^{i}_R(M, J)=0$ for some integer $i > 0$ implies
$\Ext^{j}_R(M, J)=0$ for all $j > 0$. Gluing the short exact sequence
$0\rightarrow M\rightarrow D\rightarrow M\rightarrow 0$, we get that
$M$ is  strongly $G_C$-projective.$\hfill{\square}$

{\bf Remark 3.5} From Theorem 3.3 and Proposition 3.4, we know that
every projective and $C$-projective module are strongly
$G_C$-projective. Indeed, suppose that $M$ is projective or
$C$-projective, then it is clear that $M\in \Add_R ( C\oplus R )$.
Moreover, we have a split exact sequence $0\rightarrow M\rightarrow
M\oplus M\rightarrow M\rightarrow 0$.

\vspace{0.5cm}

\centerline{\large\bf 4. Global dimensions of a ring relative to a semidualizing module}

\vspace{0.2cm}

The $\mathcal{P}_{C}$-projective and $\mathcal{I}_{C}$-injective dimension of a ring $R$ are defined as
$$\P_C\-\PD(R) = \sup\{\mathcal{P}_{C}\-\pd_R(M)| M\ is\ an\ R\-module\}$$
$$\I_C\-\ID(R) = \sup\{\mathcal{I}_{C}\-\id_R(M)| M \ is \ an \ R\-module\}.$$

When $C = R$, they are the classical homological dimensions of the ring $R$. It is natural to ask
whether the $\mathcal{P}_{C}$-projective and $\mathcal{I}_{C}$-injective dimension of a ring $R$ are equal.

{\bf Proposition 4.1} {\it For a ring $R$, $\P_{C}\-\PD(R)  = \I_{C}\-\ID(R) = \gldim(R).$}

{\it Proof.} Let $M$ be an $R$-module, then $\mathcal{P}_{C}\-\pd_R(M) = \pd_R(\Hom_R(C, M))$
from [12, Theorem 2.11]. So  $\P_{C}\-\PD(R) \leq \gldim(R)$, and the converse holds true since
$\pd_R(M) = \mathcal{P}_{C}\-\pd_R(C\otimes_R M)$. Therefore $\P_{C}\-\PD(R) = \gldim(R)$.

Similarly, we get $\I_{C}\-\ID(R) = \gldim(R)$. Thus
$\P_{C}\-\PD(R)  = \I_{C}\-\ID(R).$ $\hfill{\square}$

Then one can get a new characterization of semisimple rings in terms of $C$-projective and
$C$-injective modules.

{\bf Corollary 4.2} {\it  For a ring $R$, the following are equivalent:

(1) $R$ is semisimple,

(2) Every $R$-module is $C$-projective,

(3) Every $R$-module is $C$-injective.}

The $G_{C}$-projective and $G_{C}$-injective dimension of a ring $R$ are defined as
$$\G_{C}\-\PD(R) = \sup\{G_{C}\-\pd_R(M)| M\ is\ an\ R\-module\}$$
$$\G_{C}\-\ID(R) = \sup\{G_{C}\-\id_R(M)| M \ is \ an \ R\-module\}.$$

The following lemma is useful in the proof of the main result, and its proof
uses so-called Bass class techniques. Recall from [14] that, the Bass class with
respect to $C$, denoted $\mathcal{B}_C(R)$, consists of
all $R$-modules $N$ satisfying

(a) $\Ext_R^{i\geq1} (C, N) = 0,$

(b) Tor$^R_{i\geq1} (C, \Hom_R(C,N)) = 0,$ and

(c) The evaluation map $C\otimes_R \Hom_R(C,N)\rightarrow N$ is an isomorphism.

{\bf Lemma 4.3} {\it Let $M$ be an $R$-module with $\id_R(M)
<\infty$ and $G_C\-\pd_R(M) <\infty$. Then
$\mathcal{P}_{C}\-\pd_R(M) = G_{C}\-\pd_R(M) <\infty$.}

{\it Proof.} Since $G_C\-\pd_R(M)$ is finite, by [10, Lemma 2.8],
there is a short exact sequence of $R$-modules
$$0\rightarrow M\rightarrow N\rightarrow G\rightarrow 0\eqno{(\ast)}$$
such that $\mathcal{P}_{C}\-\pd_R(N)<\infty$ and $G$ is
$G_C$-projective.

We claim that $\Ext_R^1 (G,M) = 0$. In fact, as $G$ is $G_C$-projective,
there exists an exact sequence
$$\mathbf{X} = 0\rightarrow G\rightarrow C\otimes_RP^0\rightarrow
C\otimes_RP^1\rightarrow \cdots$$ with $P^i$ projective. Set $G_i =
\Ker(C\otimes_RP^i\rightarrow C\otimes_RP^{i+1})$ for each $i\geq
0$. The finiteness of $\id_R(M)$ implies that $M\in
\mathcal{B}_C(R)$ by [9, Corollary 6.6], and so $\Ext_R^{i\geq1}
(C,M) = 0$. Thus $\Ext_R^{i\geq1} (C\otimes_RP,M)$ $\cong$
$\Hom_R(P, \Ext_R^{i\geq1} (C,M)) = 0$ for each projective
$R$-module $P$ from [11, P.258, 9.20]. Applying $\Hom_R(-, M)$ to the sequence
$\mathbf{X}$, we have $\Ext_R^1 (G_0,M)$ $\cong$ $\Ext_R^{d+1}
(G_d,M) = 0$ by dimension-shifting argument, where $d = \id_R(M)$, as
claimed.

This implies the sequence ($\ast$) splits, then
$\sup\{\mathcal{P}_{C}\-\pd_R(M), \mathcal{P}_{C}\-\pd_R(G)\}
 = \mathcal{P}_{C}\-\pd_R(N) < \infty$, and hence $\mathcal{P}_{C}\-\pd_R(M) <\infty$.
The equality $\mathcal{P}_{C}\-\pd_R(M) = G_{C}\-\pd_R(M)$ now
follows from the result [14, Proposition 2.16].$\hfill{\square}$

{\bf Theorem 4.4} {\it For a ring $R$, $\G_{C}\-\PD(R)  =
\G_{C}\-\ID(R).$}

{\it Proof.} Assume that $\G_{C}\-\PD(R)$ is finite and not more than $n$ for some integer $n$.

Firstly, suppose that $M$ is a strongly $G_C$-projective $R$-module. We claim
that $G_C\-\id_R(M)\leq n$. By Proposition 3.4, there exists an
exact sequence $0\rightarrow M\rightarrow D\rightarrow M\rightarrow
0$, with $D\in \Add_R ( C\oplus R )$. Let $0\rightarrow M\rightarrow I_0\rightarrow I_1\rightarrow
\cdots$ be an injective resolution of $M$. By the dual version of [11,
Lemma 6.20], we have a commutative diagram
$$\begin{array}{ccccccccc}
    &   & 0 &   & 0 &   & 0 &  &   \\
     &   & \downarrow &   & \downarrow  &   & \downarrow &   &  \\
   0 & \rightarrow & M & \rightarrow & D& \rightarrow&M &\rightarrow & 0 \\
     &   & \downarrow &   & \downarrow  &   & \downarrow &   &  \\
   0 & \rightarrow & I_0 & \rightarrow & I_0\oplus I_0 & \rightarrow& I_0 & \rightarrow &0 \\
     &   & \downarrow &   & \downarrow  &   & \downarrow &   &  \\
     &   & \vdots &   & \vdots  &   & \vdots &   &  \\
     &   & \downarrow &   & \downarrow  &   & \downarrow &   &  \\
   0 & \rightarrow & I_{n-1} & \rightarrow & I_{n-1}\oplus I_{n-1} & \rightarrow& I_{n-1} & \rightarrow &0 \\
   &   & \downarrow &   & \downarrow  &   & \downarrow &   &  \\
  0 & \rightarrow & K_n & \rightarrow & G& \rightarrow&K_n &\rightarrow & 0 \\
     &   & \downarrow &   & \downarrow  &   & \downarrow &   &  \\
     &   & 0 &   & 0 &   & 0 &  &   \\
  \end{array}
$$
Because $\G_{C}\-\PD(R)\leq n$, by Lemma 2.7,
$\id_R(C\otimes_R P)\leq n$ for any projective module $P$.
It follows from [12, Theorem 2.11] that $\mathcal{I}_C\-\id_R(Q) =
\id_R(C\otimes_RQ)\leq n$ for any projective module $Q$. Thus
$G_{C}\-\id_R(C\otimes_R P\oplus Q) = \sup\{G_{C}\-\id_R(C\otimes_R P),
G_{C}\-\id_R(Q)\} $ $\leq n$ since every injective and
$C$-injective module are both $G_{C}$-injective. Therefore
$G_{C}\-\id_R(D)\leq n$, and so $G$ is $G_{C}$-injective by the dual
version of [14, Proposition 2.12].

So we obtain an exact sequence

$$\mathbf{Y}= \cdots\stackrel{f}{\rightarrow} G\stackrel{f}{\rightarrow}
G\stackrel{f}{\rightarrow} G\stackrel{f}{\rightarrow} \cdots$$ with
$G$ is $G_{C}$-injective and $K_n\cong \Ker\, f$. Applying $\Hom_R
(\Hom_R(C,I), -)$ to the sequence $\mathbf{Y}$ for any injective
module $I$, we get $\Ext^{i}_R(\Hom_R(C,I), K_n)$ $\cong$
$\Ext^{i+1}_R(\Hom_R(C,I), K_n)$ for each $i>0$. Because $I$ is
injective and $G_{C}\-\pd_R(I)\leq n$, Lemma 4.3 implies
$\mathcal{P}_{C}\-\pd_R(I)\leq n$. Then $\pd_R(\Hom_R(C,I)) =
\mathcal{P}_{C}\-\pd_R(I)\leq n$ by [12, Theorem 2.11], and So
$\Ext^{1}_R(\Hom_R(C,I), K_n)$ $\cong$ $\Ext^{n+1}_R(\Hom_R(C,I), K_n)=0$.
Thus $\Hom_R (\mathcal{I}_{C}, -)$ leaves the sequence $\mathbf{Y}$
exact, and hence $K_n$ is $G_{C}$-injective by the injective version
of Lemma 2.5. So $G_C\-\id_R(M)\leq n$ as claimed. This yields, from
[14, Proposition 2.11] and Theorem 3.3, that $G_C\-\id_R(N)\leq n$
for any $G_C$-projective $R$-module $N$.

Finally, let $M$ be an arbitrary $R$-module. By hypothesis
$G_{C}\-\pd_R(M)\leq n$. We may assume that $G_C\-\pd_R(M)\not =0$.
Then, there exists an exact sequence $0\rightarrow K \rightarrow N
\rightarrow M \rightarrow 0$ such that $N$ is $G_C$-projective and
$\mathcal{P}_{C}\-\pd_R(K)\leq n-1$ from the proof of [14, Theorem
3.6]. By induction, $G_C\-\id_R(K)\leq n$. It follows from the dual
version of [10, Lemma 3.2] that $G_C\-\id_R(M)\leq n$ since
$G_C\-\id_R(N)\leq n$.

Therefore, the right of the equality is not more than the left one, and the converse has a dual proof.
$\hfill{\square}$

In the special case $C = R$, this recovers the main result of [3, Theorem 1.1]:

{\bf Corollary 4.5} {\it $\sup\{Gpd_R(M)| M\ is\ an\ R\-module\} =
\sup\{Gid_R(M)| M \ is \ an \ R\-module\}$}.

We call the common value of the quantities in the theorem the $C$-Gorenstein
global dimension of $R$, and denote it by $\G_C\-\gldim(R)$.
Similarly, we set
$$\G_C\-\w\gldim(R) = \sup\{G_C\-\fd_R(M) |M\ is\ an\ R\-module\}.$$

{\bf Corollary 4.6} {\it The following inequalities hold:

(1) $\G_C\-\w\gldim(R) \leq \G_C\-\gldim(R),$

(2) $\G_C\-\gldim(R) \leq \gldim(R),$
and the equality holds if $\gldim(R)<\infty.$}

{\it Proof.} (1) We may assume that $\G_C\-\gldim(R) \leq n$. We
claim that every $G_C$-projective $R$-module is $G_C$-flat.
Similarly to the proof of [7, Proposition 3.4], by [14, Proposition
2.12], it suffices to show that the character module, $N^+ =
\Hom_Z(N,Q/Z)$, of every $C$-injective $R$-module $N$ has finite
$C$-projective dimension. Indeed, by the dual version of Lemma 2.7,
$\pd_R(N)\leq n$, and so $\fd_R(N)\leq n$. Then $\id_R(N^+)\leq n$
from [11, Theorem 3.52]. It follows from Lemma 4.3 that
$\mathcal{P}_C\-\pd_R(N^+)\leq n$ as desired. Therefore,
$\G_C\-\w\gldim(R) \leq \G_C\-\gldim(R)$.

(2) The inequality holds true since every projective module is $G_C$-projective.
If $\gldim(R)$ $<\infty$, then by Proposition 4.1, $\P_{C}\-\PD(R) = \gldim(R)<\infty$, the assertion follows
from [14, Proposition 2.16].$\hfill{\square}$

Now, we study the behavior of modules over rings of finite $C$-Gorenstein global dimension.

{\bf Proposition 4.7} {\it Let $R$ be a ring, and $M$ be an
$R$-module. If $\G_C\-\gldim(R) < \infty$, then

(1) $\mathcal{P}_{C}\-\pd_R(M) <\infty$ $\Leftrightarrow$ $\id_R(M) <\infty$.

(2) $\mathcal{I}_{C}\-\id_R(M) <\infty$ $\Leftrightarrow$ $\pd_R(M) <\infty$.}

{\it Proof.} (1) Simply combine Lemma 2.7 with Lemma 4.3.

(2) has a dual proof.$\hfill{\square}$

{\bf Remark 4.8} Takahashi and White [12] posed the following
question: When $R$ is a local Cohen- Macaulay ring admitting a
dualizing module and $C$ is a semidualizing $R$-module, if $M$ is an
$R$-module of finite depth such that $\mathcal{P}_{C}\-\pd_R(M)$ and
$\mathcal{I}_{C}\-\id_R(M)$ are finite, must $R$ be Gorenstein? From
Proposition 4.7 and the classical result, we know that for the rings
of finite $C$-Gorenstein global dimension, the question is positive.

\vspace{0.5cm}

{\bf Acknowledgements} The research was supported by National
Natural Science Foundation of China (Grant No. 11126092) and Science
Research Foundation of Hangzhou Dianzi University (Grant No.
KYS075610050).

\end{document}